\documentclass{article}

\newtheorem{theorem}{Theorem}[section]
\newtheorem{corollary}[theorem]{Corollary}
\newtheorem{conjecture}[theorem]{Conjecture}
\newtheorem{lemma}[theorem]{Lemma}

 \evensidemargin\oddsidemargin
\begin{document}

\title{A Family of Well-Covered Graphs with Unimodal Independence Polynomials}
\author{Vadim E. Levit and Eugen Mandrescu \\
Department of Computer Science\\
Holon Academic Institute of Technology\\
52 Golomb Str., P.O. Box 305\\
Holon 58102, ISRAEL\\
\{levitv, eugen\_m\}@hait.ac.il}
\date{}
\maketitle

\begin{abstract}
If $s_{k}$ denotes the number of stable sets of cardinality $k$ in graph $G$%
, and $\alpha (G)$ is the size of a maximum stable set, then $%
I(G;x)=\sum\limits_{k=0}^{\alpha (G)}s_{k}x^{k}$ is the \textit{independence
polynomial} of $G$ (I. Gutman and F. Harary, 1983).

J. I. Brown, K. Dilcher and R. J. Nowakowski (2000) conjectured that the
independence polynomial of a \textit{well-covered} graph $G$ (i.e., a graph
whose all maximal independent sets are of the same size) is \textit{unimodal}%
, that is, there exists some $k$ such that
\[
s_{0}\leq s_{1}\leq ...\leq s_{k-1}\leq s_{k}\geq s_{k+1}\geq ...\geq
s_{\alpha (G)}.
\]

T. S. Michael and N. Traves (2002) provided examples of well-covered graphs
whose independence polynomials are not unimodal.

A. Finbow, B. Hartnell and R. J. Nowakowski (1993) proved that under certain
conditions, any well-covered graph equals $G^{*}$ for some $G$, where $G^{*}$
is the graph obtained from $G$ by appending a single pendant edge to each
vertex of $G$.

Y. Alavi, P. J. Malde, A. J. Schwenk and P. Erd\"{o}s (1987) asked whether
for trees the independence polynomial is unimodal. V. E. Levit and E.
Mandrescu (2002) validated the unimodality of the independence polynomials
of some well-covered trees (e.g., $P_{n}^{*},K_{1,n}^{*}$, where $P_{n}$ is
the path on $n$ vertices and $K_{1,n}$ is the $n$-star graph).

In this paper we show that for any graph $G$ with $\alpha (G)\leq 4$, the
independence polynomial of $G^{*}$ is unimodal.\newline

\textbf{keywords:}\textit{\ stable set, independence polynomial, unimodal
sequence, well-covered graph.}
\end{abstract}

\section{Introduction}

Throughout this paper $G=(V,E)$ is a simple (i.e., a finite, undirected,
loopless and without multiple edges) graph with vertex set $V=V(G)$ and edge
set $E=E(G).$ The \textit{neighborhood} of a vertex $v\in V$ is the set $%
N(v)=\{w:w\in V$ \ \textit{and} $vw\in E\}$. A vertex $v$ is \textit{pendant}
if its neighborhood contains only one vertex; an edge $e=uv$ is \textit{%
pendant} if one of its endpoints is a pendant vertex. $%
K_{n},P_{n},C_{n},K_{n_{1},n_{2},...,n_{p}}$ denote respectively, the
complete graph on $n\geq 1$ vertices, the chordless path on $n\geq 1$
vertices, the chordless cycle on $n\geq 3$ vertices, and the complete $p$%
-partite graph on $n_{1}+n_{2}+...+n_{p}$ vertices.

If $G_{1},G_{2}$ are disjoint graphs, then their \textit{Zykov sum} (\cite
{Zykov}, \cite{Zykov1}) is the graph $G_{1}+G_{2}$ with $V(G_{1})\cup
V(G_{2})$ as a vertex set and
\[
E(G_{1})\cup E(G_{2})\cup \{v_{1}v_{2}:v_{1}\in V(G_{1}),v_{2}\in V(G_{2})\}
\]
as an edge set. By $G_{1}\sqcup G_{2}$ we denote the disjoint union of the
graphs $G_{1},G_{2}$, while by $\sqcup nG$ we mean the disjoint union of $%
n\geq 2$ copies of $G$.

A \textit{stable} set in $G$ is a set of pairwise non-adjacent vertices. A
stable set of maximum size will be referred to as a \textit{maximum stable
set} of $G$, and the \textit{stability number }of $G$, denoted by $\alpha
(G) $, is the cardinality of a maximum stable set in $G$. Let $s_{k}$ be the
number of stable sets in $G$ of cardinality $k$. The polynomial
\[
I(G;x)=\sum\limits_{k=0}^{\alpha (G)}s_{k}x^{k}
\]
is called the \textit{independence polynomial} of $G$, (Gutman and Harary,
\cite{GutHar}).

A finite sequence of real numbers $\{a_{0},a_{1},a_{2},...,a_{n}\}$ is said
to be \textit{unimodal} if there is some $k\in \{0,1,...,n\}$, called the
\textit{mode} of the sequence, such that
\[
a_{0}\leq a_{1}\leq ...\leq a_{k-1}\leq a_{k}\geq a_{k+1}\geq ...\geq a_{n}.
\]

As a well-known example of a unimodal sequence, we recall the following.

\begin{lemma}
\label{lem1}The sequence of binomial coefficients is unimodal, namely, for $%
n=2m$
\[
{n \choose 0}\leq ...\leq {n \choose m-1}\leq {n \choose m}\geq {n \choose m+1}%
\geq ...\geq {n \choose n},
\]
and for $n=2m+1$%
\[
{n \choose 0}\leq ...\leq {n \choose m-1}\leq {n \choose m}={n
\choose m+1}\geq {n \choose m-1}\geq ...\geq {n \choose n}.
\]
\end{lemma}

A polynomial $P(x)=a_{0}+a_{1}x+a_{2}x^{2}+...+a_{n}x^{n}$ is called \textit{%
unimodal} if its sequence of coefficients is unimodal.

\begin{figure}[h]
\setlength{\unitlength}{1cm}%
\begin{picture}(5,2)\thicklines

  \put(3.5,0){\circle*{0.29}}
  \multiput(3.5,1)(1,0){2}{\circle*{0.29}}
  \put(3.5,2){\circle*{0.29}}
  \put(3.5,0){\line(0,1){2}}
  \put(3.5,1){\line(1,0){1}}
  \put(2.5,1){\makebox(0,0){$K_{1,3}$}}

  \put(6.5,0){\circle*{0.39}}
  \put(7.5,0){\circle*{0.29}}
  \put(6.5,1){\circle*{0.39}}
  \put(6.5,2){\circle*{0.39}}
  \put(9.5,1){\circle*{0.39}}
  \multiput(7.5,1)(1,0){2}{\circle*{0.29}}
  \put(7.5,0){\circle*{0.29}}
  \put(7.5,2){\circle*{0.29}}
  \put(6.5,0){\line(1,0){1}}
  \put(6.5,1){\line(1,0){3}}
  \put(6.5,2){\line(1,0){1}}
  \put(7.5,0){\line(0,1){2}}
  \put(5.7,1){\makebox(0,0){$K_{1,3}^{*}$}}

 \end{picture}
\caption{The graphs $K_{1,3}$ and $K_{1,3}^{*}$.}
\label{fig57}
\end{figure}
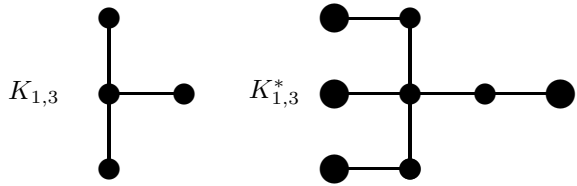

If $I(G;x)$ is unimodal, then by $mode(G)$ we mean the mode of $I(G;x)$. For
instance, the independence polynomial of $K_{1,3}$ (see Figure \ref{fig57}),
namely, $I(K_{1,3};x)=1+\mathbf{4}x+3x^{2}+x^{3}$, is unimodal and $%
mode(K_{1,3})=1$.

Hamidoune \cite{Hamidoune} proved that the independence polynomial of a
claw-free graph (i.e., a graph having no $K_{1,3}$ as an induced subgraph)
is unimodal. However, there exist graphs whose independence polynomials are
not unimodal, e.g., the graph $G=K_{100}+(\sqcup 3K_{6})$ has
\[
I(G;x)=1+\mathbf{118}x+108x^{2}+\mathbf{206}x^{3}
\]
(for other examples, see \cite{AlMalSchErdos}). Moreover, Alavi, Malde,
Schwenk and Erd\"{o}s proved in \cite{AlMalSchErdos} that for any
permutation $\sigma $ of $\{1,2,...,\alpha \}$ there is a graph $G$ with $%
\alpha (G)=\alpha $ such that $s_{\sigma (1)}<s_{\sigma (2)}<...<s_{\sigma
(\alpha )}$. Nevertheless, for trees the situation is quite different.

\begin{conjecture}
\cite{AlMalSchErdos}\label{conj1} Independence polynomials of trees are
unimodal.
\end{conjecture}

A graph $G$ is called \textit{well-covered} if all its maximal stable sets
have the same cardinality, (Plummer, \cite{Plum}). If, in addition, $G$ has
no isolated vertices and its order equals $2\alpha (G)$, then $G$ is \textit{%
very well-covered} (Favaron, \cite{Fav1}).

For $G=(V,E),V=\{v_{i}:1\leq i\leq n\}$, let $G^{*}$ be the graph obtained
from $G$ by appending a single pendant edge to each vertex of $G$, \cite
{DuttonChanBrigham}, i.e.,
\[
G^{*}=(V\cup \{u_{i}:1\leq i\leq n\},E\cup \{u_{i}v_{i}:1\leq i\leq n\}).
\]
In \cite{ToppLutz}, $G^{*}$ is denoted by $G\circ K_{1}$ and is defined as
the \textit{corona} of $G$ and $K_{1}$.

Let us notice that $G^{*}$ is well-covered (see, for instance, \cite{LevMan0}%
), and\emph{\ als}o $\alpha (G^{*})=n$. In fact, $G^{*}$ is very
well-covered, since it is well-covered, it has no isolated vertices, and its
order equals $2\alpha (G^{*})$. Moreover, the following result shows that,
under certain conditions, any well-covered graph equals $G^{*}$ for some
graph $G$.

\begin{theorem}
\cite{FinHarNow}\label{th3} Let $G$ be a connected graph of girth $\geq 6$,
which is isomorphic to neither $C_{7}$ nor $K_{1}$. Then $G$ is well-covered
if and only if its pendant edges form a perfect matching.
\end{theorem}

In other words, Theorem \ref{th3} shows that apart from $K_{1}$ and $C_{7}$,
connected well-covered graphs of girth $\geq 6$ are very well-covered. For
example, a tree $T\neq K_{1}$ could be only very well-covered, and this is
the case if and only if $T=G^{*}$ for some tree $G$ (see also Ravindra, \cite
{Ravindra}). There is a closed relationship between the independence
polynomials of $G^{*}$ and its skeleton $G$, emphasized by the following
result.

\begin{theorem}
\label{th1}\cite{LevMan4} Let $G$ be a graph of order $n$ and
\[
I(G;x)=\sum\limits_{k=0}^{\alpha (G)}s_{k}x^{k},\
I(G^{*};x)=\sum\limits_{k=0}^{\alpha (G^{*})}t_{k}x^{k}
\]
be the independence polynomials of $G$ and $G^{*}$, respectively. Then the
formulae connecting the coefficients of $I(G;x)$ and $I(G^{*};x)$ are
\begin{eqnarray*}
t_{k} &=&\sum\limits_{j=0}^{k}s_{j}\cdot {n-j \choose k-j},k\in
\{0,1,...,\alpha (G^{*})=n\}, \\
s_{k} &=&\sum\limits_{j=0}^{k}(-1)^{k+j}\cdot t_{j}\cdot {n-j
\choose n-k}
\end{eqnarray*}

In other words, $A\cdot \overline{s}=\overline{t}$, where
\begin{eqnarray*}
A &=&\left[ a_{kj}\right] ,a_{kj}={n-j \choose k-j},k\in
\{0,1,...,\alpha
(G^{*})=n\},j\in \{0,1,...,\alpha (G)\}, \\
\overline{s} &=&\left[ s_{0},s_{1},...,s_{\alpha (G)}\right] ,\overline{t}%
=\left[ t_{0},t_{1},...,t_{n}\right] .
\end{eqnarray*}
\end{theorem}

In \cite{BrownDilNow} it was conjectured that the independence polynomial of
any well-covered graph is unimodal. Recently, Michael and Traves \cite
{MichaelTraves} showed that this conjecture was true for well-covered graphs
with $\alpha (G)\in \{1,2,3\}$, and provided counterexamples for $\alpha
(G)\in \{4,5,6,7\}$. For instance, the independence polynomial
\[
1+6844x+\mathbf{10806}x^{2}+10804x^{3}+\mathbf{11701}x^{4}
\]
of the well-covered graph $\sqcup 4K_{10}+K_{1701\times 4}$ is not unimodal
(by $K_{1701\times 4}$ we mean the complete $1701$-partite graph with each
part consisting of $\sqcup 4K_{1}$). They also proposed a new conjecture,
the so-called ''\textit{roller-coaster}'' conjecture, which asserts that the
numbers
\[
s_{\left\lceil \alpha /2\right\rceil },s_{\left\lceil \alpha /2\right\rceil
+1},...,s_{\alpha (G)}
\]
of a well-covered graph $G$ with $\alpha (G)=\alpha $ are unconstrained in
the sense of Alavi \textit{et al}., and they verified their conjecture for
well-covered graphs $G$ having $\alpha (G)\leq 7$.

It is still unknown wether the independence polynomial of every very
well-covered graph is unimodal.

Using the fact that the independence polynomial of a claw-free graph is
unimodal (see Hamidoune, \cite{Hamidoune}), we verified the unimodality of
independence polynomials of a number of well-covered trees, e.g., $%
P_{n}^{*},K_{1,n}^{*},n\geq 1$ (see \cite{LevitMan2}, \cite{LevitMan3}).
These findings seem promising for proving Conjecture \ref{conj1} in the case
of very well-covered trees, because any such tree can be recursively
obtained by edge-joining of a number of $P_{n}^{*},K_{1,n}^{*}$ (see \cite
{LevitMan1}).
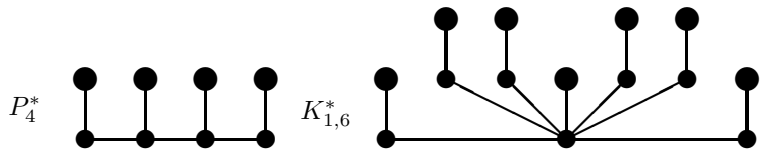
\begin{figure}[h]
\setlength{\unitlength}{0.8cm}
\begin{picture}(5,3)\thicklines

 \multiput(2,0.5)(1,0){4}{\circle*{0.29}}
  \multiput(2,1.5)(1,0){4}{\circle*{0.40}}
  \put(2,0.5){\line(0,1){1}}
  \put(3,0.5){\line(0,1){1}}
  \put(4,0.5){\line(0,1){1}}
  \put(5,0.5){\line(0,1){1}}
  \put(2,0.5){\line(1,0){3}}
  \put(1,1){\makebox(0,0){$P_{4}^{*}$}}

  \multiput(7,0.5)(3,0){3}{\circle*{0.29}}
  \multiput(7,1.5)(6,0){2}{\circle*{0.40}}
  \multiput(8,1.5)(1,0){2}{\circle*{0.29}}
  \multiput(11,1.5)(1,0){2}{\circle*{0.29}}
  \put(10,1.5){\circle*{0.40}}
  \multiput(8,2.5)(1,0){2}{\circle*{0.40}}
  \multiput(11,2.5)(1,0){2}{\circle*{0.40}}
  \multiput(7,0.5)(1,0){6}{\line(1,0){1}}
  \multiput(7,0.5)(3,0){3}{\line(0,1){1}}
  \multiput(8,1.5)(1,0){2}{\line(0,1){1}}
  \multiput(11,1.5)(1,0){2}{\line(0,1){1}}
  \put(10,0.5){\line(-2,1){2}}
  \put(10,0.5){\line(-1,1){1}}
  \put(10,0.5){\line(1,1){1}}
  \put(10,0.5){\line(2,1){2}}
  \put(6,0.9){\makebox(0,0){$K_{1,6}^{*}$}}

 \end{picture}
\caption{Well-covered trees with unimodal independence polynomials.}
\label{fig99}
\end{figure}

In this paper we show that the independence polynomials of a number of very
well-covered graphs are unimodal. More precisely, we prove that $I(G^{*};x)$
is unimodal for any $G^{*}$ whose skeleton $G$ has $\alpha (G)\leq 4$.

\section{Results}

\begin{theorem}
\label{th2}If $G$ is a graph of order $n$ and $\alpha (G)\leq 3$, then $%
I(G^{*};x)$ is unimodal with
\[
\left\lfloor \frac{n+1}{2}\right\rfloor \leq mode(G^{*})\leq \left\lfloor
\frac{n+1}{2}\right\rfloor +1.
\]

In particular, if $\alpha (G)=2$ and $n$ is odd, or $\alpha (G)=1$, then
\[
mode(G^{*})=\left\lfloor \frac{n+1}{2}\right\rfloor .
\]
\end{theorem}

\setlength {\parindent}{0.0cm}\textbf{Proof.} Let $G$ be a graph with $%
\alpha (G)=3$, and
\[
I(G;x)=s_{0}+s_{1}x+s_{2}x^{2}+s_{3}x^{3}=1+nx+s_{2}x^{2}+s_{3}x^{3}.
\]
Then $\alpha (G^{*})=n$ and $I(G^{*};x)\
=\sum\limits_{k=0}^{n}t_{k}x^{k}$, where the sequence
$t_{0},t_{1},...,t_{n}$ is given explicitly by Theorem \ref{th1}
as follows: $A\cdot \overline{s}=\overline{t}$. Let us notice
that, according to Lemma \ref{lem1}, each column of the matrix
$A=\left[ {n-j \choose k-j}\right] _{k,j}$ is a unimodal sequence.
Further in the proof
we emphasize the greatest column numbers in bold.%
\setlength
{\parindent}{3.45ex}

We analyze separately the two following cases depending on the parity of $n$.

\begin{itemize}
\item  \textit{Case 1.} $n$ is even, say $n=2m$.
\end{itemize}

We show that $t_{0}\leq t_{1}\leq ...\leq t_{m}\ $and$\ t_{m+1}\geq
t_{m+2}\geq ...\geq t_{n}$.\newline

Since
\[
\left[
\begin{array}{cccc}
{2m \choose 0} & 0 & 0 & 0 \\
&  &  &  \\
{2m \choose 1} & {2m-1 \choose 0} & 0 & 0 \\
&  &  &  \\
{2m \choose 2} & {2m-1 \choose 1} & {2m-2 \choose 0} & 0 \\
&  &  &  \\
{2m \choose 3} & {2m-1 \choose 2} & {2m-2 \choose 1} & {2m-3 \choose 0} \\
&  &  &  \\
\mathbf{\cdot } & \mathbf{\cdot } & \mathbf{\cdot } & \mathbf{\cdot } \\
&  &  &  \\
{2m \choose i} & {2m-1 \choose i-1} & {2m-2 \choose i-2} & {2m-3 \choose i-3} \\
&  &  &  \\
\mathbf{\cdot } & \cdot & \cdot & \cdot \\
&  &  &  \\
{\mathbf{2m} \choose \mathbf{m}} & {\mathbf{2m-1} \choose
\mathbf{m-1}} &
{2m-2 \choose m-2} & {2m-3 \choose m-3} \\
&  &  &  \\
{2m \choose m+1} & {\mathbf{2m-1} \choose \mathbf{m}} & {\mathbf{2m-2} \choose \mathbf{m-1}}
& {\mathbf{2m-3} \choose \mathbf{m-2}} \\
&  &  &  \\
{2m \choose m+2} & {2m-1 \choose m+1} & {2m-2 \choose m} & {\mathbf{2m-3} \choose \mathbf{m-1
}} \\
&  &  &  \\
\cdot & \cdot & \cdot & \cdot \\
&  &  &  \\
{2m \choose 2m} & {2m-1 \choose m-1} & {2m-2 \choose m-2} & {2m-3
\choose 2m-3}
\end{array}
\right] \cdot \left[
\begin{array}{c}
s_{0} \\
\\
s_{1} \\
\\
s_{2} \\
\\
s_{3}
\end{array}
\right] =\left[
\begin{array}{c}
t_{0} \\
\\
t_{1} \\
\\
t_{2} \\
\\
t_{3} \\
\\
\cdot \\
\\
t_{i} \\
\\
\cdot \\
\\
\mathbf{t}_{m} \\
\\
\mathbf{t}_{m+1} \\
\\
\mathbf{t}_{m+2} \\
\\
\cdot \\
\\
t_{2m}
\end{array}
\right] ,
\]
we get that
\begin{eqnarray*}
t_{i} &=&{n \choose i}s_{0}+{n-1 \choose i-1}s_{1}+{n-2 \choose i-2}s_{2}+%
{n-3 \choose i-3}s_{3} \\
&\leq &{n \choose i+1}s_{0}+{n-1 \choose i}s_{1}+{n-2 \choose
i-1}s_{2}+{n-3 \choose i-2}s_{3}=t_{i+1}
\end{eqnarray*}
is true for any $i\leq m-1$, and hence, $t_{0}\leq t_{1}\leq ...\leq t_{m}$.

Similarly, we infer that $t_{m+1}\geq t_{m+2}\geq ...\geq t_{2m}$.

Therefore, the sequence $\{t_{0},t_{1},...,t_{2m}\}$ is unimodal with the
only possible places for its mode $m$ or $m+1$.

\begin{itemize}
\item  \textit{Case 2.} $n$ is odd, say $n=2m+1$.
\end{itemize}

We show that $t_{0}\leq t_{1}\leq ...\leq t_{m+1}\ $and$\ t_{m+2}\geq
t_{m+3}\geq ...\geq t_{n}$.\newline

Since
\[
\left[
\begin{array}{cccc}
{2m+1 \choose 0} & 0 & 0 & 0 \\
&  &  &  \\
{2m+1 \choose 1} & {2m \choose 0} & 0 & 0 \\
&  &  &  \\
{2m+1 \choose 2} & {2m \choose 1} & {2m-1 \choose 0} & 0 \\
&  &  &  \\
{2m+1 \choose 3} & {2m \choose 2} & {2m-1 \choose 1} & {2m-2 \choose 0} \\
&  &  &  \\
\mathbf{\cdot } & \mathbf{\cdot } & \mathbf{\cdot } & \mathbf{\cdot } \\
&  &  &  \\
{2m+1 \choose i} & {2m \choose i-1} & {2m-1 \choose i-2} & {2m-2 \choose i-3} \\
&  &  &  \\
\cdot & \cdot & \cdot & \cdot \\
&  &  &  \\
{\mathbf{2m+1} \choose \mathbf{m}} & {2m \choose m-1} & {2m-1
\choose m-2} & {2m-2 \choose m-3} \\
&  &  &  \\
{\mathbf{2m+1} \choose \mathbf{m+1}} & {\mathbf{2m} \choose
\mathbf{m}} &
{\mathbf{2m-1} \choose \mathbf{m-1}} & {2m-2 \choose m-2} \\
&  &  &  \\
{2m+1 \choose m+2} & {2m \choose m+1} & {\mathbf{2m-1} \choose
\mathbf{m}} &
{\mathbf{2m-2} \choose \mathbf{m-1}} \\
&  &  &  \\
\cdot & \cdot & \cdot & \cdot \\
&  &  &  \\
{2m+1 \choose 2m+1} & {2m \choose 2m} & {2m-1 \choose 2m-1} &
{2m-2 \choose 2m-2}
\end{array}
\right] \cdot \left[
\begin{array}{c}
s_{0} \\
\\
s_{1} \\
\\
s_{2} \\
\\
s_{3}
\end{array}
\right] =\left[
\begin{array}{c}
t_{0} \\
\\
t_{1} \\
\\
t_{2} \\
\\
t_{3} \\
\\
\cdot \\
\\
t_{i} \\
\\
\cdot \\
\\
\mathbf{t}_{m} \\
\\
\mathbf{t}_{m+1} \\
\\
\mathbf{t}_{m+2} \\
\\
\cdot \\
\\
t_{2m+1}
\end{array}
\right] ,
\]
we obtain that
\begin{eqnarray*}
t_{i} &=&{n \choose i}s_{0}+{n-1 \choose i-1}s_{1}+{n-2 \choose i-2}s_{2}+%
{n-3 \choose i-3}s_{3} \\
&\leq &{n \choose i+1}s_{0}+{n-1 \choose i}s_{1}+{n-2 \choose
i-1}s_{2}+{n-3 \choose i-2}s_{3}=t_{i+1}
\end{eqnarray*}
is true for any $i\leq m$, and hence, $t_{0}\leq t_{1}\leq ...\leq t_{m+1}$.

Analogously, we deduce that $t_{m+2}\geq ...\geq t_{2m+1}$.

Therefore, the sequence $\{t_{0},t_{1},...,t_{2m+1}\}$ is unimodal with the
only possible places for its mode $m+1$ or $m+2$.

Let us observe that for $\alpha (G)=1$ and $\alpha (G)=2$, the matrix $A$
has either $2$ or $3$ columns, respectively, and the claim on the location
of the mode of $I(G^{*};x)$ follows easily. \rule{2mm}{2mm}\newline

We give now several examples covering all the possible locations of $%
mode(G^{*})$, according to Theorem \ref{th2}.

\begin{itemize}
\item  \textit{Case 1.} $\alpha (G)=1$.

$I(K_{3}^{*}{};x)=1+6x+\mathbf{9}x^{2}+4x^{3}$ has $mode(K_{3}^{*})=\left%
\lfloor (n+1)/2\right\rfloor =2$.

\item  \textit{Case 2.} $\alpha (G)=2$ and $n$ is even.

$I(P_{4}^{*}{};x)=1+8x+21x^{2}+\mathbf{22}x^{3}+8x^{4}$ and $%
mode(P_{4}^{*})=\left\lfloor (n+1)/2\right\rfloor +1=3$, while $%
I((K_{4}-e)^{*};x)=1+8x+\mathbf{19}x^{2}+18x^{3}+6x^{4}$ has $%
mode((K_{4}-e)^{*})=\left\lfloor (n+1)/2\right\rfloor =2$.

\item  \textit{Case 3.} $\alpha (G)=2$ and $n$ is odd.

$I(P_{3}^{*}{};x)=1+6x+\mathbf{10}x^{2}+5x^{3}$ has $mode(P_{3}^{*})=\left%
\lfloor (n+1)/2\right\rfloor =2$.

\item  \textit{Case 4.} $\alpha (G)=3$.

$I(\sqcup 3K_{1}^{*};x)=1+6x+\mathbf{12}x^{2}+8x^{3}$ has $mode(\sqcup
3K_{1}^{*})=\left\lfloor (n+1)/2\right\rfloor =2$. $I((K_{1}\sqcup
P_{3})^{*};x)=1+8x+22x^{2}+\mathbf{25}x^{3}+10x^{4}$ has $mode((K_{1}\sqcup
P_{3})^{*})=\left\lfloor (n+1)/2\right\rfloor +1=3$.
\end{itemize}

The next result extends Theorem \ref{th2} for skeletons with the stability
number equal to $4$.

\begin{theorem}
\label{th4}If $G$ is a graph of order $n$ and $\alpha (G)=4$, then $%
I(G^{*};x)$ is unimodal with
\[
\left\lfloor \frac{n+1}{2}\right\rfloor \leq mode(G^{*})\leq \left\lfloor
\frac{n+1}{2}\right\rfloor +2.
\]

Moreover, if $n$ is odd, then
\[
\left\lfloor \frac{n+1}{2}\right\rfloor \leq mode(G^{*})\leq \left\lfloor
\frac{n+1}{2}\right\rfloor +1.
\]
\end{theorem}

\setlength {\parindent}{0.0cm}\textbf{Proof.} Let $G$ be a graph with $%
\alpha (G)=4$, and
\[
I(G;x)=s_{0}+s_{1}x+s_{2}x^{2}+s_{3}x^{3}+s_{4}x^{4}=1+nx+s_{2}x^{2}+s_{3}x^{3}+s_{4}x^{4}.
\]
Then $\alpha (G^{*})=n$ and $I(G^{*};x)\ =\sum\limits_{k=0}^{n}t_{k}x^{k}$,
where the sequence $t_{0},t_{1},...,t_{n}$ is given explicitly by Theorem
\ref{th1} as follows: $A\cdot \overline{s}=\overline{t}$.%
\setlength {\parindent}{3.45ex} Let us notice that, according to
Lemma \ref{lem1}, each column of the matrix $A=\left[ {n-j \choose
k-j}\right] _{k,j}$ is a unimodal sequence. Further in the proof
we emphasize the greatest column numbers in bold.

We distinguish between the two following cases depending on the parity of $n$%
.

\begin{itemize}
\item  \textit{Case 1.} $n$ is odd, say $n=2m+1$.
\end{itemize}

We show that $t_{0}\leq t_{1}\leq ...\leq t_{m+1}\ $and$\ t_{m+2}\geq
t_{m+3}\geq ...\geq t_{n}.$

Since
\[
\left[
\begin{array}{ccccc}
{2m+1 \choose 0} & 0 & 0 & 0 & 0 \\
&  &  &  &  \\
{2m+1 \choose 1} & {2m \choose 0} & 0 & 0 & 0 \\
&  &  &  &  \\
{2m+1 \choose 2} & {2m \choose 1} & {2m-1 \choose 0} & 0 & 0 \\
&  &  &  &  \\
{2m+1 \choose 3} & {2m \choose 2} & {2m-1 \choose 1} & {2m-2 \choose 0} & 0 \\
&  &  &  &  \\
{2m+1 \choose 4} & {2m \choose 3} & {2m-1 \choose 2} & {2m-2 \choose 1} & {2m-3 \choose 0} \\
&  &  &  &  \\
\mathbf{\cdot } & \mathbf{\cdot } & \mathbf{\cdot } & \mathbf{\cdot } &  \\
&  &  &  &  \\
{2m+1 \choose i} & {2m \choose i-1} & {2m-1 \choose i-2} & {2m-2
\choose i-3} & {2m-3 \choose i-4} \\
&  &  &  &  \\
\cdot & \cdot & \cdot & \cdot &  \\
&  &  &  &  \\
{\mathbf{2m+1} \choose \mathbf{m}} & {2m \choose m-1} & {2m-1
\choose m-2} & {2m-2 \choose m-3} & {2m-3 \choose m-4} \\
&  &  &  &  \\
{\mathbf{2m+1} \choose \mathbf{m+1}} & {\mathbf{2m} \choose
\mathbf{m}} & {\mathbf{2m-1} \choose \mathbf{m-1}} & {2m-2 \choose
m-2} & {2m-3 \choose m-3}
\\
&  &  &  &  \\
{2m+1 \choose m+2} & {2m \choose m+1} & {\mathbf{2m-1} \choose
\mathbf{m}} &
{\mathbf{2m-2} \choose \mathbf{m-1}} & {\mathbf{2m-3} \choose \mathbf{m-2}} \\
&  &  &  &  \\
{2m+1 \choose m+3} & {2m \choose m+2} & {2m-1 \choose m+1} & {2m-2
\choose m} &
{\mathbf{2m-3} \choose \mathbf{m-1}} \\
&  &  &  &  \\
\cdot & \cdot & \cdot & \cdot & \cdot \\
&  &  &  &  \\
{2m+1 \choose 2m+1} & {2m \choose 2m} & {2m-1 \choose 2m-1} &
{2m-2 \choose 2m-2} & {2m-3 \choose 2m-3}
\end{array}
\right] \cdot \left[
\begin{array}{c}
s_{0} \\
\\
s_{1} \\
\\
s_{2} \\
\\
s_{3} \\
\\
s_{4}
\end{array}
\right] =\left[
\begin{array}{c}
t_{0} \\
\\
t_{1} \\
\\
t_{2} \\
\\
t_{3} \\
\\
t_{4} \\
\\
\cdot \\
\\
t_{i} \\
\\
\cdot \\
\\
\mathbf{t}_{m} \\
\\
\mathbf{t}_{m+1} \\
\\
\mathbf{t}_{m+2} \\
\\
\mathbf{t}_{m+3} \\
\\
\cdot \\
\\
t_{2m+1}
\end{array}
\right] ,
\]
we get that
\[
t_{i}={n \choose i}s_{0}+{n-1 \choose i-1}s_{1}+{n-2 \choose
i-2}s_{2}+{n-3 \choose i-3}s_{3}+{n-4 \choose i-4}s_{4}\leq
\]
\smallskip
\[
\leq {n \choose i+1}s_{0}+{n-1 \choose i}s_{1}+{n-2 \choose
i-1}s_{2}+{n-3 \choose i-2}s_{3}+{n-4 \choose i-3}s_{4}=t_{i+1}
\]
is true for any $i\leq m$, i.e., $t_{0}\leq t_{1}\leq ...\leq t_{m+1}$.

Analogously, we infer that $t_{m+2}\geq ...\geq t_{2m+1}$.

Therefore, the sequence $\{t_{0},t_{1},...,t_{2m+1}\}$ is unimodal with the
only possible places for its mode $m+1$ or $m+2$.

\begin{itemize}
\item  \textit{Case 2.} $n$ is even, say $n=2m$.
\end{itemize}

We show that
\[
t_{0}\leq t_{1}\leq ...\leq t_{m},t_{m+2}\geq ...\geq t_{n},\ and\
2t_{m+1}\geq t_{m}+t_{m+2}.
\]
Hence, it follows that $t_{m+1}\geq \min \{t_{m},t_{m+2}\}$, and
consequently, $I(G^{*};x)$ is unimodal with the mode $j\in \{m,m+1,m+2\}$.

Since
\[
\left[
\begin{array}{ccccc}
{2m \choose 0} & 0 & 0 & 0 & 0 \\
&  &  &  &  \\
{2m \choose 1} & {2m-1 \choose 0} & 0 & 0 & 0 \\
&  &  &  &  \\
{2m \choose 2} & {2m-1 \choose 1} & {2m-2 \choose 0} & 0 & 0 \\
&  &  &  &  \\
{2m \choose 3} & {2m-1 \choose 2} & {2m-2 \choose 1} & {2m-3 \choose 0} & 0 \\
&  &  &  &  \\
{2m \choose 4} & {2m-1 \choose 3} & {2m-2 \choose 2} & {2m-3 \choose 1} & {2m-4 \choose 0} \\
&  &  &  &  \\
\mathbf{\cdot } & \mathbf{\cdot } & \mathbf{\cdot } & \mathbf{\cdot } &
\mathbf{\cdot } \\
&  &  &  &  \\
{2m \choose i} & {2m-1 \choose i-1} & {2m-2 \choose i-2} & {2m-3
\choose i-3} &
{2m-4 \choose i-4} \\
&  &  &  &  \\
\cdot & \cdot & \cdot & \cdot & \cdot \\
&  &  &  &  \\
{\mathbf{2m} \choose \mathbf{m}} & {\mathbf{2m-1} \choose
\mathbf{m-1}} &
{2m-2 \choose m-2} & {2m-3 \choose m-3} & {2m-4 \choose m-4} \\
&  &  &  &  \\
{2m \choose m+1} & {\mathbf{2m-1} \choose \mathbf{m}} & {\mathbf{2m-2} \choose \mathbf{m-1}}
& {\mathbf{2m-3} \choose \mathbf{m-2}} & {2m-4 \choose m-3} \\
&  &  &  &  \\
{2m \choose m+2} & {2m-1 \choose m+1} & {2m-2 \choose m} & {\mathbf{2m-3} \choose \mathbf{m-1}}
 & {\mathbf{2m-4} \choose \mathbf{m-2}} \\
&  &  &  &  \\
\cdot & \cdot & \cdot & \cdot & \cdot \\
&  &  &  &  \\
{2m \choose 2m} & {2m-1 \choose 2m-1} & {2m-2 \choose 2m-2} &
{2m-3 \choose 2m-3} & {2m-4 \choose 2m-4}
\end{array}
\right] \cdot \left[
\begin{array}{c}
s_{0} \\
\\
s_{1} \\
\\
s_{2} \\
\\
s_{3} \\
\\
s_{4}
\end{array}
\right] =\left[
\begin{array}{c}
t_{0} \\
\\
t_{1} \\
\\
t_{2} \\
\\
t_{3} \\
\\
t_{4} \\
\\
\cdot \\
\\
t_{i} \\
\\
\cdot \\
\\
\mathbf{t}_{m} \\
\\
\mathbf{t}_{m+1} \\
\\
\mathbf{t}_{m+2} \\
\\
\cdot \\
\\
t_{2m}
\end{array}
\right] ,
\]
we get that
\[
t_{i}={n \choose i}s_{0}+{n-1 \choose i-1}s_{1}+{n-2 \choose
i-2}s_{2}+{n-3 \choose i-3}s_{3}+{n-4 \choose i-4}s_{4}\leq
\]
\smallskip
\[
\leq {n \choose i+1}s_{0}+{n-1 \choose i}s_{1}+{n-2 \choose
i-1}s_{2}+{n-3 \choose i-2}s_{3}+{n-4 \choose i-3}=t_{i+1}
\]
is true for any $i+1\leq m$, and $t_{0}\leq t_{1}\leq ...\leq t_{m}$.

Similarly, we infer that $t_{m+2}\geq t_{m+3}\geq ...\geq t_{2m}$.

Now,
\[
2t_{m+1}-t_{m}-t_{m+2}=
\]
\begin{eqnarray*}
&=&\sum\limits_{i=0}^{4}s_{i}\cdot \left\{ 2{2m-i \choose
m-i+1}-{2m-i \choose m-i}-{2m-i \choose m-i+2}\right\} = \\
&=&\sum\limits_{i=0}^{4}s_{i}\cdot \frac{(2m-i)!}{m!\cdot (m-i+2)!}\cdot
(2m+3i-i^{2}-2).
\end{eqnarray*}

Hence, we see that for $0\leq i\leq 3$ the coefficients near $s_{i}$ are
non-negative. If $i=4$ and $m\geq 3$, then the coefficient near $s_{4}$ is
non-negative, too. Therefore, for $m\geq 3$ either $t_{m+1}\geq t_{m}$ or $%
t_{m+1}\geq t_{m+2}$, and th\emph{ese facts ensure th}at the sequence $%
\{t_{0},t_{1},...,t_{2m}\}$ is unimodal and the only possible places for its
mode are $m,m+1$ and $m+2$.

If $m\leq 2$, then $n\in \left\{ 2,4\right\} $.\ Since $\alpha \left(
G\right) =4$, we are left with only one graph under consideration, namely, $%
G=\sqcup 4K_{1}$. Now the theorem follows from the unimodality of the
polynomial
\[
I((\sqcup 4K_{1})^{*};x)=(1+2x)^{4}=1+8x+24x^{2}+\mathbf{32}x^{3}+16x^{4},
\]
and the fact that $mode((\sqcup 4K_{1})^{*})=3$. \rule{2mm}{2mm}\newline

We give now several examples covering all the possible locations of $%
mode(G^{*})$, according to in Theorem \ref{th4}.

\begin{itemize}
\item  \textit{Case 1.} $n$ is odd.

$I(K_{1,4}^{*};x)=$ $1+10x+36x^{2}+\mathbf{62}x^{3}+52x^{4}+17x^{5},$

$mode(K_{1,4}^{*})=\left\lfloor (n+1)/2\right\rfloor =3$;

$I(\sqcup 3K_{2}^{*}\sqcup K_{1}^{*};x)=1+14x+81x^{2}+250x^{3}+443x^{4}+%
\mathbf{450}x^{5}+243x^{6}+54x^{7},mode(\sqcup 3K_{2}^{*}\sqcup
K_{1}^{*})=\left\lfloor (n+1)/2\right\rfloor +1=5$.

\item  \textit{Case 2.} $n$ is even.

$I(K_{2}^{*};x)=$ $1+\mathbf{4}x+3x^{2},mode(K_{2}^{*})=\left\lfloor
(n+1)/2\right\rfloor =1$.

$I((\sqcup 4K_{1})^{*};x)=(1+2x)^{4}=$ $1+8x+24x^{2}+\mathbf{32}%
x^{3}+16x^{4},$

$mode(\sqcup 4K_{1}^{*})=\left\lfloor (n+1)/2\right\rfloor +1=3$.

$I(\sqcup
2P_{8}^{*};x)=1+32x+466x^{2}+4100x^{3}+24405x^{4}+104292x^{5}+331314x^{6}+799092x^{7}+1480780x^{8}+2118240x^{9}+%
\mathbf{2334666}x^{10}+$

$%
1964532x^{11}+1238901x^{12}+566780x^{13}+177610x^{14}+34100x^{15}+3025x^{16},
$

$mode(\sqcup 2P_{8}^{*})=\left\lfloor (n+1)/2\right\rfloor +2=10$.
\end{itemize}

Let us notice that the unimodality of $I(G^{*};x)$ does not imply the
unimodality of $I(G;x)$, even for well-covered graphs. For example, the
independence polynomial of $G=(\sqcup 4K_{10})+K_{1800\times 4}$ is not
unimodal:
\begin{eqnarray*}
I(G;x) &=&1800\cdot (1+x)^{4}+(1+10x)^{4}-1800 \\
&=&1+7240x+\mathbf{11400}x^{2}+11200x^{3}+\mathbf{11800}x^{4}.
\end{eqnarray*}
However, according to Theorem \ref{th4}, $I(G^{*};x)$ is unimodal, since $%
\alpha (G)=4$.

Clearly, for any $n\geq 1$, there exists a disconnected very well-covered
graph $G$ with $\alpha (G)=n$, whose $I(G;x)$ is unimodal, e.g., $G=\sqcup
nK_{2}$ has $I(\sqcup nK_{2};x)=(1+2x)^{n}$. Moreover, for any $n\geq 1$,
there exists a very well-covered tree, namely $P_{n}^{*}$, whose
independence polynomial is unimodal (see \cite{LevitMan2}). The following
result shows that this assertion is also true for a series of connected very
well-covered graphs that are not\emph{\ }trees.

\begin{corollary}
For any $n\geq 3$, there exists a connected very well-covered graph $H$
different from a tree, such that $\alpha (H)=n$ and $I(H;x)$ is unimodal.
\end{corollary}

\setlength {\parindent}{0.0cm}\textbf{Proof.} Let $1\leq m\leq \min
\{4,n-2\},G=(\sqcup mK_{1})+K_{n-m}$, and $H=G^{*}$. Then $G$ is a connected
graph different from a tree, and, consequently, $H$ is a connected non-tree,
too. Moreover, $H$ is very well-covered, and $\alpha (H)=n$. Since $\alpha
(G)=m\leq 4$, the independence polynomial $I(H;x)$ is unimodal, according to
Theorems \ref{th2} and \ref{th4}. \rule{2mm}{2mm}%
\setlength
{\parindent}{3.45ex}

\section{Conclusions}

In this paper we showed that the independence polynomial of a
number of very well-covered graphs is unimodal. This gives support
for the still open conjecture concerning the unimodality of
independence polynomials of very well-covered graphs. We leave as
an open problem the question whether $%
I(G^{*};x)$ is unimodal whenever $\alpha (G)\geq 5$.

\end{document}